\documentclass[12pt]{amsart}
\usepackage{amsmath,amscd,amssymb,amsfonts}
\newcommand{\h}{\hbox}
\newcommand{\q}{\quad}

\newcommand{\nin}{\par\noindent}
\newcommand{\bs}{\par\bigskip}
\newcommand{\ms}{\par\medskip}
\newcommand{\sk}{\par\smallskip}
\newcommand{\mtim}{\h{$\times$}}
\newcommand{\C}{{\mathbf C}}
\newcommand{\N}{{\mathbf N}}
\newcommand{\HH}{{\mathbf H}}
\newcommand{\PP}{{\mathbf P}}
\newcommand{\Z}{{\mathbf Z}}
\newcommand{\QQ}{\overline{\mathbf Q}}
\newcommand{\Hc}{{\mathcal H}}

\newcommand{\f}{{}\,\overline{\!f}{}}

\newcommand{\pio}{\overline{\pi}{}}
\newcommand{\X}{{}\,\overline{\!X}{}}
\newcommand{\Y}{\overline{Y}{}}
\newcommand{\ZZ}{\overline{Z}}
\newcommand{\W}{\overline{W}}

\newcommand{\YY}{\widetilde{Y^{\rm qf}}}
\newcommand{\pit}{\widetilde{\pi}}
\newcommand{\eo}{\overline{\eta}}
\newcommand{\Spec}{{\rm Spec}}
\newcommand{\Ext}{{\rm Ext}}
\newcommand{\MHS}{{\rm MHS}}
\newcommand{\VMHS}{{\rm VMHS}}
\newcommand{\NFad}{{\rm NF}^{\rm ad}}
\newcommand{\sm}{{\rm sm}}
\newcommand{\bl}{\bigl}
\newcommand{\br}{\bigr}
\newcommand{\into}{\hookrightarrow}
\newcommand{\simto}{\buildrel{\sim\,\,}\over\longrightarrow}
\newcommand{\ssc}{\,\raise.15ex\h{${\scriptstyle\circ}$}\,}

\newcommand{\les}{\leqslant}
\begin{document}
\title{Normal functions and spread of zero locus}
\author{Morihiko Saito}
\address{RIMS Kyoto University, Kyoto 606-8502 Japan}
\begin{abstract}
If there is a topologically locally constant family of smooth algebraic varieties together with an admissible normal function on the total space, then the latter is constant on any fiber if this holds on some fiber. Combined with spreading out, it implies for instance that an irreducible component of the zero locus of an admissible normal function is defined over k if it has a k-rational point where k is an algebraically closed subfield of the complex number field with finite transcendence degree. This generalizes a result of F. Charles that was shown in case the normal function is associated with an algebraic cycle defined over k.
\end{abstract}
\maketitle
\centerline{\bf Introduction}
\bs\nin
Let $k$ be an algebraically closed subfield of $\C$ with finite transcendence degree. Let $X$ be a smooth complex variety defined over $k$. Let $\HH$ be an admissible variation of mixed Hodge structure of strictly negative weights on $X$ (see \cite{Ka}, \cite{SZ}). Let $\nu$ be an admissible normal function of $\HH$, which is a holomorphic section of the family of Jacobians $J(\HH)$ satisfying some good properties (see \cite{GGK}, \cite{adm}). Let $Z$ be an irreducible component of the zero locus $\nu^{-1}(0)$ of $\nu$. This is algebraic as a corollary of \cite{BP}, \cite{KNU}, \cite{Sch2} (see \cite{BPS}).
\sk
Assume $Z$ is not defined over $k$. Let $K$ be the (minimal) field of definition of $Z$. This is the smallest subfield $K\subset\C$ containing $k$ and such that $Z$ is defined over $K$, see \cite[Cor.~4.8.11]{Gro2} (and also \cite{We}). Let $R$ be a finitely generated $k$-subalgebra of $K$ whose field of fractions is $K$. Set
$$S:=\Spec\,R\otimes_k\C.$$
The dimension of $S$ coincides with the relative transcendence degree of $K$ over $k$, and is called the {\it transcendence degree} of $Z$ over $k$. We may assume that $S$ is sufficiently small by replacing $R$ without changing $K$. Then $S$ is smooth, and there is a closed subvariety
$$Y\subset X\times S,$$
defined over $k$ and having the morphisms
$$f:Y\to S,\q\pi:Y\to X,$$
induced by the projections and such that $\pi$ induces an isomorphism
$$Y_{s_0}:=f^{-1}(s_0)\simto Z\subset X,$$
where $s_0$ is the $k$-generic point of $S$ corresponding to the inclusion $R\into K\into\C$. (In this paper a point of a complex algebraic variety means always a closed point.)
\sk
Let $\ZZ^k$ denote the $k$-Zariski closure of $Z$ in $X$, i.e. the smallest Zariski closed subset of $X$ defined over $k$ and containing $Z$. This coincides with the Zariski closure of $\pi(Y)\subset X$, see Remark~(2.2)(iii) below. Set $Y_s:=Y\cap\bl(X\times\{s\}\br)\subset X$ for $s\in S$. We have the following.
\ms\nin
{\bf Theorem~1.} {\it Assume $S$ is sufficiently small. Then the induced morphism $\pi:Y\to\ZZ^k$ is dominant and quasi-finite so that $\dim\ZZ^k=\dim Z+\dim S$. Moreover there is a Zariski closed subset $\Sigma\subset S$ such that $\{s_0\}$ is a connected component of $\Sigma$ and we have $Y_s\subset\nu^{-1}(0)$ if $s\in\Sigma$, and $Y_s\subset X\setminus\nu^{-1}(0)$ if $s\notin\Sigma$.}
\ms
Here no condition is assumed about the relation between the normal function $\nu$ and the field $k$, although there is a rather strong restriction coming from the assumption that the non-constant normal function $\nu$ is defined on a smooth complex algebraic variety $X$ defined over $k$ (see Remark~(3.8)(i) below).
The quasi-finiteness of $\pi$ implies that the restriction of the family over any curve $C$ on $S$ has no fixed point (locally on $C$). It might be possible that $Y_s$ for some $s\in\Sigma\setminus\{s_0\}$ is contained in another irreducible component of $\nu^{-1}(0)$ intersecting $Z$ (although no such examples are explicitly known, see also Remark~(3.8)(ii) below). Note that $\Sigma$ is smooth by Theorem~2 below (by shrinking $S$ if necessary). However, it is unclear whether every irreducible component of the zero locus has the same dimension.
\sk
For the proof of Theorem~1, we show the following.
\ms\nin
{\bf Proposition~1.} {\it In the notation and the assumption of Theorem~$1$, there is a resolution of singularities $Y'\to Y$ such that the induced morphism $f':Y'\to S$ is smooth and topologically locally trivial over $S$, and the pull-back of the normal function $\nu$ to $Y'$ by $\pi'$ coincides with the pull-back of some normal function of $f'_*\pi'{}^*\HH$ by $f'$, where $f',\pi'$ denote the compositions of $Y'\to Y$ with $f,\pi$ respectively.}
\ms
Here $f'_*\pi'{}^*\HH$ is an admissible variation of mixed Hodge structure (see \cite{mhm}). This assertion easily follows from the theory of admissible normal function explained in Section 1.
\sk
Theorem~1 implies the following corollary (which was shown by F.~Charles \cite{Ch2} using a completely different method \cite{Ch1} if $\nu$ is associated with an algebraic cycle defined over $k$).
\ms\nin
{\bf Corollary~1.} If $Z$ contains a $k$-rational point $z$, then $Z$ is defined over $k$.
\ms
In the case $Z$ is a point (which is important for the proof of the Hodge conjecture \cite{Sch1}) Corollary~1 and Theorem~1 say nothing. In fact, $S$ is a $k$-Zariski open subset of $\ZZ^k$, and $\pi$ is the natural inclusion in this case.
\sk
In the situation of Theorem~1, we say that $Z$ has {\it maximal transcendence degree} over $k$ if $\dim S={\rm codim}_XZ$. Note that the last condition is equivalent to $X=\ZZ^k$ if $X$ is irreducible.
\ms\nin
{\bf Theorem~2.} {\it In the notation of Theorem~$1$, assume $Z$ has maximal transcendence degree over $k$. Then $Z$ is smooth, and does not intersect the other irreducible components of $\nu^{-1}(0)$.}
\ms
I would like to thank F.~Charles, B.~Kahn, and M.~Levine for useful discussions about the subject of this paper. In fact, the original idea was inspired by the discussions with them.
This work is partially supported by Kakenhi 24540039.
\sk
In Section 1 we review some basics of admissible normal functions.
In Section 2 we recall the notion of spread out.
In Section 3 we prove the main theorems after recalling some facts from algebraic or analytic geometry.
\bs\bs
\vbox{\centerline{\bf 1. Admissible normal functions}
\bs\nin
In this section we review some basics of admissible normal functions.}
\ms\nin
{\bf 1.1.~Absolute case.} Let $X$ be a smooth complex algebraic variety. Let $\HH$ be an admissible variation of mixed Hodge structure of strictly negative weights on $X$, see \cite{Ka}, \cite{SZ}. (Here the underlying $\Z$-local system is assumed torsion-free.) We have the family of Jacobians $J(\HH)$ over $X$. Its fiber at $x\in X$ is set-theoretically identified with
$$J(\HH_x)\,\bl(:=\HH_{x,\C}/(F^0\HH_{x,\C}+\HH_{x,\Z})\br)=\Ext^1_{\MHS}(\Z,\HH_x),$$
where $\HH_x$ is the fiber of $\HH$ at $x$, and MHS is the abelian category of graded-polarizable mixed $\Z$-Hodge structures \cite{De} (see \cite{Ca} for the last isomorphism).
\sk
Let $\NFad(X,\HH)$ be the group of admissible normal functions of $\HH$ (see \cite{adm} and also \cite{GGK}). We have a canonical isomorphism
$$\NFad(X,\HH)=\Ext^1_{\VMHS(X)}(\Z_X,\HH),$$
where $\VMHS(X)$ is the abelian category of admissible variations of mixed Hodge structure on $X$ (assumed always graded-polarizable). We have the following exact sequence (see e.g. \cite[Theorem~3.6]{ext}):
$$0\to J\bl(H^0(X,\HH)\br)\to\NFad(X,\HH)\buildrel{cl\,}\over\to H^1(X,\HH),
\leqno(1.1.1)$$
where $H^i(X,\HH)\in\MHS$ by \cite{mhm}. It denotes also its underlying $\Z$-module as in the case of the last term of (1.1.1). The image of $\nu$ by $cl$ in (1.1.1) is called the cohomology class of $\nu$, and is defined by considering the extension class of underlying local systems.
\sk
If $cl(\nu)=0$, then $\nu$ is the pull-back of an element of the first term of (1.1.1) by $X\to pt$, and is called {\it constant} on $X$.
\ms\nin
{\bf 1.2.~Relative case.} Let $f:Y\to S$ be a smooth morphism of smooth connected complex algebraic varieties. Set $Y_s:=f^{-1}(s)$ for $s\in S$. Let $\HH$ be an admissible variation of mixed Hodge structure of strictly negative weights on $Y$. Assume $f$ is topologically locally trivial over $S$. Then $\Hc^0f_*\HH$ is an admissible variation of mixed Hodge structure of strictly negative weights on $S$ with
$$f^*\Hc^0f_*\HH\subset\HH,$$
and $\Hc^1f_*\HH$ is an admissible variation of mixed Hodge structure on $S$ (by using the stability of mixed Hodge modules by cohomological direct images, see \cite{mhm}). Moreover, there is an exact sequence
$$0\to\NFad(S,\Hc^0f_*\HH)\buildrel{f^*}\over\to\NFad(Y,\HH)\to\Gamma(S,\Hc^1f_*\HH),
\leqno(1.2.1)$$
where the last term means the global sections of the underlying local system of $\Hc^1f_*\HH$.
The image of an admissible normal function $\nu$ of $\HH$ to the last term is the cohomology class of $\nu$, and is denoted by $cl(\nu)$. The exactness of (1.2.1) is proved by taking the restriction to each fiber of $f$ and reducing to (1.1.1). If $cl(\nu)$ vanishes, then $\nu$ is an admissible normal function of $f^*\Hc^0f_*\HH\subset\HH$, and coincides with the pull-back of some admissible normal function of $\Hc^0f_*\HH$.
\ms\nin
{\bf Proposition~1.3.} {\it Let $\nu$ be an admissible normal function of $\HH$. If the restriction of $\nu$ to the fiber $Y_{s_0}$ of $f$ at some $s_0\in S$ is constant in the sense of $(1.1)$, then the restriction of $\nu$ to any fiber $Y_s$ is constant, and $\nu$ is the pull-back of some admissible normal function of $\Hc^0f_*\HH$.}
\ms\nin
{\it Proof.} The assumption implies that the section of $\Hc^1f_*\HH$ corresponding to $cl(\nu)$ vanishes at $s_0$, and hence $cl(\nu)=0$. So the assertion follows from the exact sequence (1.2.1).
\bs\bs
\vbox{\centerline{\bf 2. Spread of closed subvarieties}
\bs\nin
In this section we recall the notion of spread out.}
\ms\nin
{\bf 2.1.~Spread out.} Let $X$ be a complex algebraic variety. Let $k$ be an algebraically closed subfield of $\C$ with finite transcendence degree. We say that $X$ is defined over $k$, if there is a $k$-variety $X_k$ endowed with an isomorphism
$$X=X_k\otimes_k\C,$$
where $\otimes_k\C$ means the base change by $\Spec\,\C\to\Spec\,k$. We say that a point $x$ of $X$ is {\it $k$-generic} if there is no proper closed subvariety $V$ of an irreducible component of $X$ such that $x\in V$ and $V$ is defined over $k$, see also \cite{We}. We say that $U$ is a {\it $k$-Zariski open} subset if its complement is a closed subvariety defined over $k$.
\sk
Let $X$ be a complex algebraic variety defined over $k$, and $Z$ be a closed complex algebraic subvariety of $X$. There is a complex smooth affine variety $S$ together with a closed subvariety $Y$ of $X\mtim S$ which are all defined over $k$ and such that $Z$ coincides with the fiber
$$Y_{s_0}:=Y\cap\bl(X\mtim\{s_0\}\br)\subset X,$$
for some $k$-generic point $s_0$ of $S$. We say that $Y$ is a {\it spread} of $Z$ over $k$. We will denote respectively by $\pi$ and $f$ the morphisms of $Y$ induced by the first and second projections from $X\mtim S$.
\sk
More precisely, $Z\subset X$ is defined over some subfield $K$ of $\C$ which is finitely generated over $k$. Let $R$ be a finitely generated $k$-subalgebra of $K$ whose filed of fractions is $K$. Set $S_k:=\Spec\,R$, $S:=\Spec\,R\otimes_k\C$. If $R$ is sufficiently large, then there is a closed $k$-subvariety
$$Y_k\into X_{S_k}:=X_k\mtim_k S_k,$$
such that its base change by $\eo:\Spec\,\C\to S_k$ is identified with the natural inclusion $Z\into X$. Here $\eo$ is induced by the inclusions $R\into K\into\C$, and induces a morphism $R\otimes_k\C\to\C$ defining $s_0\in S$. We may assume $S$ smooth by shrinking $S_k$ if necessary.
\ms\nin
{\bf Remarks~2.2.} (i) There is the smallest subfield $K_Z\subset\C$ such that $Z$ is defined over it, see \cite[Cor.~4.8.11]{Gro2}. This field $K_Z$ is called the {\it minimal field of definition} of $Z\subset X$ over $k$.
\ms
(ii) The minimum of the dimension of $S$ for all the spreads of $Z$ is called the {\it transcendence degree} of $Z$ over $k$, and is denoted by ${\rm tr.\,deg}_k\,Z$. This coincides with ${\rm tr.\,deg}_k\,K_Z$, that is the relative transcendence degree over $k$ of the minimal field of definition $K_Z$.
\ms
(iii) Assume $\dim S$ is minimal (i.e. $\dim S={\rm tr.\,deg}_k\,Z$). Then the field of fractions of the affine ring $R$ of $S_k$ is a finite extension of the minimal field of definition $K_Z$. Moreover we have the equality
$$\overline{\pi(Y)}=\ZZ^k,$$
where the left-hand side is the closure of $\pi(Y)$ in $X$, and the right-hand side is the $k$-Zariski closure of $Z$ in $X$, that is the minimal closed $k$-subvariety of $X$ containing $Z$.
Indeed, take the closure of $Y$ in $X\times\overline{S}$ with $\overline{S}$ a compactification of $S$ defined over $k$. Its image by the projection $X\times\overline{S}\to X$ is a $k$-Zariski closed subvariety containing $\pi(Y)$ as a dense subset (even in the classical topology). So the Zariski closure of $\pi(Y)$ is defined over $k$, and hence contains the $k$-Zariski closure of $Z$. On the other hand, we can replace $X$ with the $k$-Zariski closure of $Z$ for the construction of spread of $Z$. This implies the opposite inclusion. So the desired equality follows.
\ms\nin
{\bf Examples~2.3.} (i) Assume $k=\QQ$, $X=\C^2$ and $Z=\{y=cx\}\subset\C^2$ where $c$ is a transcendental number. In this case $Z$ contains the origin, which is a $k$-rational point, and the fibers $Y_s$ of the spread $Y$ of $Z$ are given by $\{y=sx\}\subset\C^2$ where $s$ is identified with the coordinate of $S=\C$.
\ms
(ii) Assume $k=\QQ$, $X=\C^3$, and $Z$ is the line $C(P_c)$ passing through the origin and the point
$$P_c:=(1-c^2,2c,1+c^2)\in\C^3.$$
Here $c$ is a transcendental number and $C(P_c)$ means the cone of $[P_c]\in{\bf P}^2$. In this case we have
$$\ZZ^k=\{x^2+y^2=z^2\}\subset\C^3,$$
and the fibers $Y_s$ of the spread $Y$ of $Z$ are given by the lines $C(P_s)$ where $s$ is identified with the coordinate of $S=\C$ and $P_s$ is defined similarly to $P_c$. This is closely related to the parametrization of the projective curve $\{x^2+y^2=z^2\}\subset\PP^2$ by $\PP^1$ using the intersection with the lines $\{y=s(x+z)\}$ for $s\in\C$.
\bs\bs
\vbox{\centerline{\bf 3. Proof of the main theorems}
\bs\nin
In this section we prove the main theorems after recalling some facts from algebraic or analytic geometry.}
\ms\nin
{\bf 3.1.~Some consequence of Weierstrass preparation theorem.} Let $\pi:Y\to X$ be a morphism of complex analytic spaces. Assume $y\in Y$ is an isolated point of $\pi^{-1}\pi(y)$. Then there is an open neighborhood $V_y$ of $y$ in $Y$ which is finite over an open neighborhood $U_x$ of $x:=\pi(y)$ in $X$.
\sk
Indeed, we may assume that $X$ is a closed analytic subset of a polydisk $\Delta^m$, and $Y$ is a closed analytic subset of $\Delta^n\mtim X$ defined by holomorphic functions $h_1,\dots,h_r$, where $\pi:Y\to X$ is induced by the projection $\Delta^{n+m}\to\Delta^m$ (by using the graph embedding), and the origins of the polydisks $\Delta^{n+m}$ and $\Delta^n$ respectively correspond to $y$ and $x$. Since $y$ is isolated in $\pi^{-1}(x)$, we may assume that the restriction of $h_1$ to $\Delta^n\mtim\{0\}$ is not identically zero, and then $h_1(x_1,0,\dots,0)$ is not identically zero. Let $W$ be the closed subset of $\Delta^{n+m}$ defined by $h_1$. By the Weierstrass preparation theorem, the projection $W\to\Delta^{n+m-1}$ is a finite morphism by shrinking the polydisks if necessary. Moreover, the image of an analytic space by a finite (or more generally, projective) morphism is a closed analytic subspace (by using the direct image of the structure sheaf). So the assertion follows by induction on $n$. (If $\pi:Y\to X$ is algebraically defined, this assertion also follows from Remark~(ii) below.)
\ms\nin
{\bf 3.2.~Grothendieck's version of Zariski's main theorem.} Let $\pi:Y\to X$ be a quasi-projective morphism of algebraic varieties defined over $k$. By \cite[Thm.~4.4.3]{Gro1} the quasi-finite locus $Y^{\rm qf}$ of $\pi$ is a $k$-Zariski open subset of $Y$. Here the quasi-finite locus means the subset consisting of $y\in Y$ which is isolated in $\pi^{-1}\pi(y)$. (Note that this assertion is local on $X$.) Moreover the morphism $\pi^{\rm qf}:Y^{\rm qf}\to X$ induced by $\pi$ can be extended to a finite morphism $\pit:\YY\to X$ defined over $k$ in such a way that $Y^{\rm qf}$ is a dense $k$-Zariski open subvariety of $\YY$. This is Grothendieck's version of Zariski's main theorem.
\ms\nin
{\bf 3.3.~Proof of Proposition~1.} Let $Y,S$ be as in (2.1) where we assume that $\dim S$ is minimal. Let $\X$ be a smooth compactification of $X$ defined over $k$. Let $\Y$ be the closure of
$$Y\subset X\mtim S\q\h{in}\q\X\mtim S,$$
with $\pio$, $\f$ the morphisms of $\Y$ induced by the first and second projections from $\X\mtim S$. These are also defined over $k$. Note that $\pio^{\,-1}(X)=Y$ since $Y$ is closed in $X\mtim S$.
\sk
Let $\rho:\Y'\to\Y$ be a resolution of singularities defined over $k$ and inducing an isomorphism over the smooth part $Y^{\sm}$ of $Y$. We may assume that the inverse image of $D:=\Y\setminus Y^{\sm}$ by $\rho$ is a divisor with simple normal crossings on $\Y'$ such that any intersections of irreducible components of $\rho^{-1}(D)$ are smooth over $S$ by assuming $S$ sufficiently small. Here we may also assume that $\rho^{-1}(\Y\setminus Y)$ is a divisor by using a blow-up. Then the above assertion holds also for $\rho^{-1}(\Y\setminus Y)$, since the latter is a divisor contained in $\rho^{-1}(\Y\setminus Y^{\sm})$.
\sk
Set $\HH':=\pi'{}^*\HH$, and $\nu':=\pi'{}^*\nu$. Then $\Hc^1f'_*\HH'$ is an admissible variation of mixed Hodge structure on $S$, and $cl(\nu')$ vanishes so that $\nu'$ comes from an admissible normal function $\nu''$ of $f'_*\HH'$, see Remark~(1.3). So the assertion follows. This finishes the proof of Proposition~1.
\ms\nin
{\bf 3.4.~Proof of Theorem~1.} Let $\Sigma$ be the zero locus of $\nu''$ on $S$. Note that $s_0\in\Sigma$ and the zero locus of $\nu'$ coincides with $f^{-1}(\Sigma)$.
So the fiber $Y_s\subset X$ of $f$ over any $s\in\Sigma$ is contained in $\nu^{-1}(0)$, and $Y_s$ does not intersect $\nu^{-1}(0)$ if $s\notin\Sigma$.
Note that this implies that $Y_{s_0}$ is contained in the quasi-finite locus $Y^{\rm qf}$ of $\pi:Y\to X$ if $\Sigma$ is $0$-dimensional at $s_0$. We can show the last property as follows.
\sk
If $\Sigma$ has positive dimension at $s_0$, then $Y_s=Z$ for any $s$ in the connected component $\Sigma_0$ of $\Sigma$ containing $s_0$ (since $Y_s$ must contain at least a non-empty smooth open subset of $Z$ in the classical topology by using Proposition~1). Take a general hyperplane section $H$ defined over $k$ and intersecting $\Sigma_0$ (by using the fact that any non-empty Zariski open subset of complex projective space has a $k$-rational point). Then we may replace $S$ with the smallest closed subvariety of $H$ defined over $k$ and containing a point of $H\cap \Sigma_0$. But this contradicts the minimality of $\dim S$. So $\Sigma$ must be $0$-dimensional at $s_0$.
\sk
It now remains to show the quasi-finiteness of $\pi:Y\to X$ by assuming $S$ sufficiently small. We have $Y_{s_0}\subset Y^{\rm qf}$ by the above argument. So the proof of Theorem~1 is reduced to Lemma~(3.5) below with $W$ the non-quasi-finite locus $Y\setminus Y^{\rm qf}$, since the latter is a closed subvariety of $Y$ defined over $k$ by (3.2). (Note that finitely generated fields are countable and algebraic closures of countable fields are also countable.)
\ms\nin
{\bf Lemma~3.5.} {\it Let $f:Y\to S$ be a morphism of complex algebraic varieties defined over a countable field $k$, where $S$ is irreducible. Let $W$ be a closed subvariety of $Y$ defined over $k$, and $s_0$ be a $k$-generic point of $S$. Assume $W\cap f^{-1}(s_0)=\emptyset$. Then $f(W)$ is contained in a proper closed subvariety of $S$ defined over $k$.}
\ms\nin
{\it Proof.} We may assume $W$ irreducible. The hypothesis implies that
$$W\cap f^{-1}(s)=\emptyset,$$
for any $k$-generic point $s\in S$ by using the action of ${\rm Aut}(\C/k)$. Let $S''$ denote the union of all the proper closed subvarieties of $S$ defined over $k$. This is a countable union, i.e.,
$$\h{$S''=\bigcup_{i\in\N} S''_i$},$$
with $S''_i$ proper closed subvarieties of $S$ defined over $k$ ($i\in\N$), since $k$ is a countable field. Moreover the set of $k$-generic points of $S$ coincides with $S\setminus S''$ (see also \cite[Prop.~3.2]{SS}). So we get
$$W\subset f^{-1}(S'').$$
Hence $W=W\cap f^{-1}(S'')$ is a countable union of closed subvarieties
$$W_i:=W\cap f^{-1}(S''_i)\subset W\q(i\in\N).$$
This implies that $W=W_i$ for some $i$, and the assertion follows. This finishes the proofs of Lemma~(3.5) and Theorem~1.
\ms\nin
{\bf 3.6.~Proof of Corollary~1.} Assume $Z$ is not defined over $k$. Apply Theorem~1 to $Z$, where $\dim S>0$. We then see that $s_0$ in Theorem~1 is the image of a $k$-rational point of $Y$ by the morphism $f:Y\to S$. For this, consider
$$(z,s_0)\in X\times S,$$
which belongs to $Y$, and is moreover $k$-rational. Indeed, it is a connected component of the intersection defined over $k$
$$Y\cap\bl(\{z\}\mtim S\br)\subset X\mtim S,$$
by using the quasi-finiteness of $\pi:Y\to X$ (and $k$ is algebraically closed). This is a contradiction since $s_0$ is a $k$-generic point of $S$. So the assertion follows.
\ms\nin
{\bf 3.7.~Proof of Theorem~2.} We have $\dim X=\dim Y$ by assumption. Let $D'\subset X$ denote the discriminant locus of $\pi':Y'\to X$ in Proposition~1, i.e. $D'$ is the smallest Zariski closed subset such that $\pi'$ is \'etale over the complement of $D'$. This is defined over $k$. Note also that $D'$ contains the image of the singular locus of $Y$. Since the $k$-Zariski closure of $Z$ is $X$, we get $Z\not\subset D'$, i.e.,
$$Z\setminus D'\ne\emptyset.
\leqno(3.7.1)$$
\sk
Assume $Z$ is singular. Let $z$ be any singular point of $Z$, and $z'$ be the corresponding point of $Y_{s_0}$ by $\pi$. In the classical topology, there is an open neighborhood $V_{z'}$ of $z'$ in $Y$ which is finite over an open neighborhood $U_z$ of $z$ in $X$ by the morphism $\pi_z:V_{z'}\to U_z$, induced by $\pi$, see (3.1). The morphism $\pi_z$ is locally biholomorphic over $U_z\setminus D'$, and induces an isomorphism over $Z\cap U_z$ by using Theorem~1. Here we assume $V_{z}$ sufficiently small so that
$$f(V_{z'})\cap\Sigma=\{s_0\}.$$
Then the degree of $\pi_z$ must be 1 since $(U_z\cap Z)\setminus D'\ne\emptyset$ by (3.7.1). So $V_{z'}$ must be smooth. Applying this argument to any singular point of $Z$, we see that $Y$ is smooth on a neighborhood of $Y_{s_0}$. However, this contradicts the assumption that the geometric generic fiber $Y_{s_0}$ of $f:Y\to S$ (which is isomorphic to $Z$) is singular. Thus $Z$ is smooth, and the first assertion is proved.
\sk
It is then easy to see that $Z$ cannot intersect the other components of $\nu^{-1}(0)$, since $\pi(Y)$ contains an open neighborhood of $Z$ (at least in the classical topology). This finishes the proof of Theorem~2.
\ms\nin
{\bf Remarks~3.8.} (i) The quasi-finiteness of $\pi$ in Theorem~1 does not hold for the spread of an arbitrary closed irreducible subvariety $Z$ of a smooth complex algebraic variety $X$ defined over $k$, see Example~(2.3)(i).
In this example, any admissible variation of mixed Hodge structure on $X=\C^2$ is constant, and so is any admissible normal function on $X$.
Here we cannot replace $X$ with the self-product of an elliptic curve $E$, since ${\rm End}(E)$ is discrete (more precisely, it is isomorphic to $\Z$ or $\Z^2$).
If we replace $X$ with the ``cone" of a $k$-Zariski open subset of $\PP^1$ having a non-constant normal function on it, then $X$ cannot contain the origin (as long as it is a variety), and we get no contradictions.
Note also that the cone of an abelian variety is quite singular although it is complete and has many non-constant normal functions on it.
\ms
(ii) There is an example of a closed subvariety $W$ of a smooth complex algebraic variety $X$ defined over $k$ such that $W\subset\W^k_1$ with $W_1$ an irreducible component of $W$ and some fiber $Y_s$ of the spread $Y$ of $W_1$ coincides with another irreducible component of $W$ intersecting $W_1$. Consider, for instance, the case $k=\QQ$, $X=\C^3$ and $W=\bigcup_{0\les i<r}W_i$ with $W_i=C(P_{c_i})$ in the notation of Example~(2.3)(ii), where the $c_i$ are transcendental numbers. In this case we have
$$\W^k_i=\{x^2+y^2=z^2\}\subset\C^3.$$
Note that there is no non-constant admissible normal function on $X$.

\end{document}